\documentclass[12pt]{article}

\usepackage{amsbsy}
\usepackage{amsmath}
\usepackage{amssymb}
\usepackage{epsfig}
\usepackage{multicol}
\usepackage{multirow}
\usepackage{color}
\usepackage{bm}
\usepackage{threeparttable}
\usepackage{float}
\renewcommand{\baselinestretch} {1.3}
\makeatletter 
\def\singlespace{\def\baselinestretch{1}\@normalsize}

\@addtoreset{equation}{section}

\renewcommand{\hat}{\widehat}

\newcommand{\qed}{$\Box$}

\newcommand{\bg}{\begin{eqnarray}}
\newcommand{\ed}{\end{eqnarray}}
\newcommand{\bgn}{\begin{eqnarray*}}
\newcommand{\edn}{\end{eqnarray*}}

\makeatletter
\def\singlespace{\def\baselinestretch{1}\@normalsize}

\date{\today}

\title{The strong representation for the nonparametric estimator of  length-biased and right-censored data}

\author{
Jianhua Shi$^{1,3}$, Xiaoping Chen$^{2,3*}$ and Yong Zhou$^{3}$\\
 {\footnotesize $^{1}$School of Mathematics and Statistics,}
          {\footnotesize Minnan Normal University, Zhangzhou, China }\\
   {\footnotesize $^{2}$School of Mathematics and Computer Science,}
          {\footnotesize Fujian Normal University, Fuzhou, China }\\
          {\footnotesize $^{3}$Academy of Mathematics and Systems Science,}
          {\footnotesize Chinese Academy of Sciences, Beijing  China, }\\
           {\footnotesize  and School of Statistics and Management,}
          {\footnotesize Shanghai University of Finance and Economics, Shanghai, China }\\
 }

\begin{document}
\maketitle

\begin{singlespace}
\begin{footnotetext}
{\textbf{Acknowledgements}: Zhou's work was supported by National Natural Science Foundation of China (NSFC) (71271128), the State Key Program of National Natural Science Foundation of China (71331006), NCMIS, Shanghai First-class Discipline A and IRTSHUFE. Shi and Chen's work was supported by Natural Science Foundation of Fujian Province, China (2012J01028)
and National Natural Science Foundation of China (NSFC) (11301473).

\textbf{*Corresponding author:}
Xiaoping Chen, Fujian Normal University, No.1 Keji Road,  Shangjie University Town, Minhou County, Fuzhou,  Fujian, China 350117.
 (\textit{E-mail: xpchen@fjnu.edu.cn}).
}
\end{footnotetext}
\end{singlespace}

\begin{abstract}
In this paper, we consider the modified product-limit estimator of an unknown distribution function proposed by Huang and Qin(2011), where the observations are subject to length-biased and right-censored data. A strong representation result for the modified product-limit estimator is  established with a remainder  $O(n^{ - {3 \mathord{\left/ {\vphantom {3 4}} \right. \kern-\nulldelimiterspace} 4}} (\log n)^{{3 \mathord{\left/  {\vphantom {3 4}} \right. \kern-\nulldelimiterspace} 4}} ) $ a.s. Such results are very useful when we consider statistics that are the functional of the estimator of nonparametric distribution function. Also, an uniform consistency rate of the estimator is given.

\end{abstract}


\baselineskip=20pt
\section{Introduction}

Length-biased data frequently appear in observational studies, when the observed samples are not randomly selected
from the population of interest but with probability proportional to their length  (Shen \emph{et al.} 2009), such as  the prevalent sampling design which only considers subjects with disease.   In the prevalent sampling, the truncation time has a uniform distribution since the occurrence of disease onset follows a stationary Poisson process. Hence, the length-biased data are associated with left-truncation data  as those who fail before  sampling time are not observable. In addition to being length-biased, survival sampling data are usually subject to right censoring due to loss of follow-up.
Various  methods for estimating the distribution function
have been developed when the distribution function of truncation time is unspecified. Lots of work has been devoted to the strong representation for the distribution function estimator under left-truncation and right-censored (LTRC) data but very few under length-biased and right-censored (LBRC) setting.

A brief review of the theoretic development for the LTRC data in recent years is worth mentioning.
Assume that $(T^0,A,C^0)$ is a random vector, where $T^0$ is the survival time of
interest with unknown cumulative distribution function (c.d.f) $F ( \cdot )$, $A$ is a random
left truncation time with unknown c.d.f. $F_A ( \cdot )$ and $C^0$ is a random
right censored time with arbitrary c.d.f. $F_C ( \cdot )$. Let $Y^0 = \min
(T^0,C^0)$, and $\Delta = I(T^0 \le C^0)$ be the usual indicator of censoring status.
Denote $\bar R(t) = n^{ - 1}\sum\nolimits_{i = 1}^n {I(a_i \le t \le y_i )} $ as the empirical estimator for $R(t) = P(\left. {A \le t \le Y^0} \right)$, where and throughout the paper, the lowercase letters for the random variables indicate the sampling value from a population.
Tsai \emph{et al.} (1987) proposed the nonparametric estimator $\hat {F}_n( \cdot) $ for $F ( \cdot)$,
\[
1 - \hat {F}_n (x) = \prod\limits_{y_i \le x} {[1 - (n\bar R (y_i ))^{ -
1}]^{\delta _i }},
\]
which is the well-known TJW product-limit (PL) estimator.
Obviously, the estimator reduces to the Lynden-Bell (1971) PL-estimator for only left truncation and to the Kaplan-Meier (1958) PL-estimator for purely
right-censored data ($T=0$).

Extensive literature focused on the strong representation for the TJW PL-estimator is available. Cs\"{o}rg\H{o} and Horv$\acute{a}$th  (1982a,1983) and Burke \emph{et al.} (1981) studied the strong representation of the Kaplan-Meier estimator for right-censored data. For general censored, Cs\"{o}rg\H{o} and Horv$\acute{a}$th  (1982b,c) established strong approximations for the Kaplan-Meier estimator.
  Under LTRC sampling, Lai and Ying  (1991) obtained a functional law of the iterated
logarithm for the modified PL-estimator by using martingale theory. Furthermore, 
Gijbels and Wang (1993) studied a strong
representation for TJW PL-estimator, and obtained the order $O(n^{ - 1}  \log n  )$ a.s. of the remainder term under some suitable conditions. Zhou (1996)
considered more complicated situation for distribution function support,
and proved a strong approximation for TJW PL process $\{\hat F_n (t) - F(t), t>0\}$  at the rate $O(n^{ - 1} \log ^{1 + \varepsilon } n)\ a.s.$
The rate of approximation is improved by Zhou and Yip (1999), where the remainder term  is of order
$O(n^{ -1}\log \log n)$ a.s. under some suitable integrability assumptions, which is currently the best result about the convergent rate. For references to some other relevant strong
approximation results of TJW PL-estimator see Sellero \emph{et al.} (2005), Liang \emph{et al.} (2009), Liu \emph{et al.} (2013) among others.

Therefore the strong representation of TJW PL-estimator is an interesting problem in the field of probability and statistics,
which attracts much attention. However, the strong behavior for the PL-estimator in LBRC
data is relatively less studied. Recently, using the potential property in the LBRC
design, Huang and Qin (2011) obtained a weak representation for a modified  TJW PL-estimator of an unknown survival distribution. The negligible term in the representation is $o_p(1)$.

In this article, the remainder term in the representation is considered. Utilizing the tool of empirical process, we obtain two almost sure representations for a cumulative hazard estimator and the estimator of the corresponding distribution function, respectively. The negligible terms in the representations are firstly verified of order $O(n^{{{ - 3} \mathord{\left/
 {\vphantom {{ - 3} 4}} \right.
 \kern-\nulldelimiterspace} 4}} (\log n)^{{3 \mathord{\left/
 {\vphantom {3 4}} \right.
 \kern-\nulldelimiterspace} 4}} ) $ almost sure.
 Such results are very useful when we consider statistics that are the functional of the estimator of nonparametric distribution function, such as investigating the properties of quantile function of the modified  TJW PL-estimator,  studying the oscillation modulus of the estimator and  estimating density  function, and ROC curve, etc. For instance, Cs\"{o}rg\H{o} and Horv$\acute{a}$th (1983) investigated the maximal deviation of the PL- estimate from the estimated distribution function, Lo \emph{et al.} (1989) studied the estimation of a density and a hazard rate function by a strong  uniform approximation of the Kaplan-Meier estimator. For more discussion and application about strong representation, one may refer to Burke \emph{et al.} (1981,1988), Cs\"{o}rg\H{o} and Horv$\acute{a}$th (1982),  Horv$\acute{a}$th (1984), Gijbels and Wang (1993) and Zhou and Yip (1999), Tse (2003), etc. In comparison, such statistical analysis is hard to be performed based on the weak representation in Huang and Qin (2011).

  The rest of the paper is organized as follows. Some
notations and the main theoretical results concerning two nonparametric estimators are presented in Section 2. Section 3 is devoted to several prepared lemmas and their proofs.

\section{Notations and Main Results} 
We now introduce related random variable notations for LBRC sampling, Huang and Qin's (2011) notation is followed whenever possible.  Let $(T^0,A,C^0)$
denote a random vector where $T^{\rm{0}}$ is the interested survival time from the
disease incidence to the failure event with marginal density function $f(t)$
and survival function $S(t)$; $A$ is the random left truncation time from
the disease incidence to sampling time, $\xi $, and $C^0$ is the total censoring time
from the disease onset. Meanwhile, suppose $W^0$ be the onset time for the disease incidence, $C$ be the time from
sampling time to censoring,\ i.e. it is the residual censored time, then $C^0=A+C$.

Two basic assumptions for the general population, needed throughout the paper, are presented as follows.

(A1) The distribution of $T^0$ is independent of $W^0$.

(A2) The incidence of disease onset occurs over calendar time at a
constant rate, that is, $W^0$ has a constant density function.

The above assumptions will not be stated again for the sake of simplicity.
For obtaining the strong approximation of the remainder term, an additional integrability hypothesis is needed. To this end, let $F^u(t) = P(\Delta = 1,Y \le t)$ denote the subdistribution function, and define $a_G  = \inf \{ t:G(t) > 0\}, b_G  = \sup \{ t:G(t) < 1\}$ for any d.f. $G(\cdot)$. Assume $Y \sim H( \cdot )$, then we have $1-H(\cdot)=(1-F(\cdot))(1-F_C(\cdot))$ and $b_H$  = min$ (b_F , b_{F_C })$ by the independence assumption. Compared with Woodroofe's (1985) results, $F(\cdot)$ can be reconstructed if $a_{F_A }  \le a_H$ and $ b_{F_A }  \le b_H$. Therefore we assume that $a_{F_A }  \le a_H$ and $ b_{F_A }  \le b_H$ shall hold throughout this paper. Meanwhile,  put $0/ 0=0$ for convenience.

(A3) For $a_H  < b < b_H$, $\int_{a_H }^b {R^{-3} (u)dF^u(u)}  < \infty .$

\textbf{\emph{Remark 1}} The Assumption (A3) is satisfied when $a_{F_A}<a_H$, which is similar to the condition (2.1) in Zhou and Yip (1999). When the random variables $A$ and $C$ are independent, (A3) reduces to the condition (0.2) of Stute (1993) with truncated data.

 Next for illustrating the left truncated sampling, we drop the superscript
$^0$ in the notation of $W^0$, $Y^0$ and $T^0$, and thus $ (W,T) \sim _d \left. {(W^0,T^0)} \right|T^0 \ge \xi - W^0 > 0,$
where $ \sim _d $ denotes identical distribution. In our setting, $C$ is assumed to be independent of $(W,\xi ,T)$ and $\xi $ is
independent with $(W,T)$. 
 Define $\alpha = :P(Y \ge A)$, just
as the notations defined before, when $Y < A$ nothing is
observed in LBRC model. Naturally one needs to assume $\alpha > 0$. Furthermore, if we set $\tilde {V} = \min (V,C)$, where $V$ is the residual survival time
from the sampling time, then the observed data is i.i.d. copies of
$(W,A,\tilde {V},\Delta )$. In this article, all the quantities with $\tilde {}$  are the new empirical processes considered.

Denote the survival functions of the random variables $A , T, C$ and $V$ defined
in the prevalent population as $S_A (t),S_T (t),S_C (t)$, and $S_V (t)$, which are always assumed continuous in the paper, and
the corresponding marginal density functions as $f_A (t),f_T (t),f_C (t)$ and $f_V (t)$, respectively. Under (A1) and (A2), an important relation is that the truncation time and the residual survival
time share the same marginal density function, i.e.
\[
f_A (t) = f_V (t) = \frac{S(t)}{\mu }I(t > 0),
\]
where  $\mu = E(T^0)= \int_0^\infty {uf(u)du} .$

Based on the key property, Huang and Qin (2011) proposed to replace the empirical estimate $\bar {R}(t)$ by the estimator
\[
\tilde {R}(t) = n^{ - 1}\sum\limits_{j = 1}^n {I(y_j \ge t)} - \tilde {S}_A
(t),\]
  where
$\tilde {S}_A (t) = \prod\limits_{u \in [0,t]} {\{1 - \tilde {K}^{-1}(u) d\tilde
{Q}(u)\}} $ is the usual Kaplan--Meier estimator for $A$,
\[
\tilde {Q}(t) = \frac{1}{n}\sum\limits_{i = 1}^n {[I(a_i \le t) + \delta _i
I(\tilde {v}_i \le t)]} = :\tilde Q_1 (t) + \tilde Q_2 (t) ,
\]
and
\[
\tilde {K}(t) = \frac{1}{n}\sum\limits_{i =
1}^n {[I(a_i \ge t) + I(\tilde {v}_i \ge t)]}{ = :\tilde K_1 (t) + \tilde K_2 (t)}.
\]

Let $\Lambda (\cdot) $ be the cumulative hazard
function of $F ( \cdot )$. Note that $E[dF^u(t)]={\mu }^{-1}f(t)\int_0^t {G(s)ds}$ under LBRC mechanism. It is easy to see that
$\Lambda (t) = \int_0^t {R^{-1}(u)dF^u(u)} .$ Hence, using the estimator $\tilde {R}(t)$ above, an alternative nonparametric estimator for $\Lambda (\cdot)$ can be constructed by combining properties from both $A$ and $V$ under length-biased sampling,
\[
\tilde {\Lambda }(t) = \int_0^t {\frac{d\bar {N}(u)}{\tilde {R}(u)}} ,
\]
where $\bar {N}(t) = n^{ - 1}\sum\limits_{j = 1}^n {\delta _j I(y_j \le t)}$, and the corresponding c.d.f. estimator is $\tilde {F}_n(t)$ with
\[1 - \tilde {F}_n(t) =  \prod\limits_{u \in [0,t]} {\{ 1 - d\tilde \Lambda (u)\} }.\]
 Huang and Qin (2011) presented a weak asymptotic
large sample property, expressing the corresponding survival estimator as i.i.d. means of random variables with a negligible term.

Denote $Q(t) = E\{ I(A \le t) + \Delta I(\tilde V \le t)\}  = :Q_1 (t) + Q_2 (t)$ and $K(t) = E\{ I(A \ge t) + I(\tilde V \ge t)\}  = :K_1 (t) + K_2 (t)$, respectively.  Now, define  i.i.d. stochastic processes for $1 \le i \le n$,
 \[
\phi _i (t) = \int_0^t {K^{ - 2}(u)\{I(a_i \ge u) + I(\tilde {v}_i \ge u)\}}
dQ(u) - \frac{I(a_i \le t)}{K(a_i )} - \frac{\delta _i I(\tilde {v}_i \le
t)}{K(\tilde {v}_i )}.
\]
It can be shown that $\{\phi _i (t), 1 \le i \le n\}$ is  a mean zero stochastic process.

 Write $\psi _i (t) = \psi _{1i} (t) + \psi _{2i} (t)$, where
\[
\psi _{1i} (t)\mbox{ = }\int_0^t {R^{ - 2}(u)I(y_i \ge u \ge a_i )dF^u(u)} -
\frac{\delta _i I(y_i \le t)}{R(y_i )},
\]
and
\[
\psi _{2i} (t)\mbox{ = }\int_0^t {R^{ - 2}(u)\{I(a_i > u) - S_A (u) - S_A
(u)\phi _i (u)\}dF^u(u)} .
\]

The strong asymptotic representations of the proposed estimator $\tilde {\Lambda }(\cdot)$ and
$\tilde F_n (\cdot)$ are summarized in Theorem 2.1 and Theorem 2.2, respectively. Assume $a_H=0$ in the proofs w.l.o.g. throughout to avoid trivialities.

{\bf {\it \textbf{Theorem 2.1}.}} Suppose that (A3) holds for some $b < b_H$. Then uniformly in $a_H \le t \le b < b_H ,$
the stochastic process $\tilde {\Lambda }(t) - \Lambda
(t)$ has an asymptotic representation
\[
\tilde \Lambda (t) - \Lambda (t) =  - n^{ - 1} \sum\limits_{i = 1}^n {[\psi _{1i} (t) + \psi _{2i} (t)]}  + {R_{n1} (t)}
\]
with $
\mathop {\sup }\limits_{a_H  \le t \le b} \left| {R_{n1} (t)} \right| = O(n^{{{ - 3} \mathord{\left/
 {\vphantom {{ - 3} 4}} \right.
 \kern-\nulldelimiterspace} 4}} (\log n)^{{3 \mathord{\left/
 {\vphantom {3 4}} \right.
 \kern-\nulldelimiterspace} 4}} )\ a.s.$

\textbf{\emph{Proof of Theorem 2.1.}} Utilizing Lemma 3.3 and by following the discussion of Lemma 3.1 we have
\begin{eqnarray*}
 &&\int_0^t {(\frac{1}{\tilde {R}(u)}} - \frac{1}{R(u)})d(\bar {N}(u) -
F^u(u))= O(n^{{{ - 3} \mathord{\left/
 {\vphantom {{ - 3} 4}} \right.
 \kern-\nulldelimiterspace} 4}} (\log n)^{{3 \mathord{\left/
 {\vphantom {3 4}} \right.
 \kern-\nulldelimiterspace} 4}} )
 \ a.s.
 \end{eqnarray*}

Meanwhile,
\begin{eqnarray*}
&& \int_0^t {\frac{{(R(u) - \tilde R(u))^2 }}{{\tilde R(u)R^2 (u)}}} dF^u (u) \le \mathop {\sup }\limits_{0 \le u \le t} (R(u) - \tilde R(u))^2 \int_0^t {\frac{dF^u (u)}{{\tilde R(u)R^2 (u)}}} 
=O(n^{ - 1} \log n)\ a.s.
\end{eqnarray*}

By the definition of $\psi _{2i} (t)$  and Remark 3 below Lemma 3.3, we have
 \begin{eqnarray*}
  - n^{ - 1} \sum\limits_{i = 1}^n {\psi _{2i} (t)}
     &=&  - \int_0^t {\frac{{\tilde R(u) - \bar R(u)}}{{R^2 (u)}}dF^u (u)}  + O(n^{ - {3 \mathord{\left/ {\vphantom {3 4}} \right. \kern-\nulldelimiterspace} 4}} log^{{3 \mathord{\left/ {\vphantom {3 4}} \right.
 \kern-\nulldelimiterspace} 4}}n )\ a.s.
 \end{eqnarray*}

Then, one can easily decompose the term $\tilde {\Lambda }(t) - \Lambda (t)$.
\begin{eqnarray*}
 \tilde {\Lambda }(t) - \Lambda (t) &=& \int_0^t {\frac{d(\bar {N}(u) - F^u(u))}{\tilde {R}(u)}} + \int_0^t
{(\frac{1}{\tilde {R}(u)}} - \frac{1}{R(u)})dF^u(u) \\
    & = & - \frac{1}{n}\sum\limits_{i = 1}^n {[\psi _{1i} (t) + \psi _{2i} (t)]}  + O(n^{{{ - 3} \mathord{\left/
 {\vphantom {{ - 3} 4}} \right.
 \kern-\nulldelimiterspace} 4}} (\log n)^{{3 \mathord{\left/
 {\vphantom {3 4}} \right.
 \kern-\nulldelimiterspace} 4}} )
\ a.s.
 \end{eqnarray*}

This completes the proof of Theorem 2.1.   \ \ \ \ \ \ \ \ \ \ \ \ \ \ \ \ \ \ \ \ \  \ \ \  \ \ \ \ \ \ \ \ \ \ \ \ \ \ \ \ \ \ \ \ \  \ \ \  \ \ \ \ \ \ \ \ \ \ \ \ \ \ \ \ \ \ \ \ \  \ \ \   \qed

\textbf{\emph{Theorem 2.2}}. Suppose that (A3) holds for some $b < b_H$, then we have uniformly in $a_H \le t \le b < b_H ,$
 \begin{eqnarray*}
 && \tilde F_n (t) - F(t)  = n^{ - 1}\sum\limits_{i = 1}^n {(1-F(t))[\psi _{1i} (t) + \psi _{2i} (t)]} +  {R_{n2} (t)}
 \end{eqnarray*}
with $
\mathop {\sup }\limits_{a_H  \le t \le b} \left| {R_{n2} (t)} \right| = O(n^{{{ - 3} \mathord{\left/
 {\vphantom {{ - 3} 4}} \right.
 \kern-\nulldelimiterspace} 4}} (\log n)^{{3 \mathord{\left/
 {\vphantom {3 4}} \right.
 \kern-\nulldelimiterspace} 4}} )\ a.s.$

\emph{\textbf{Remark 2}} The approximation rate firstly
obtained in the paper is mainly based on Lemma 3.1. We have not got a
more appropriate tool at present. The approximation rate may be improved, which is an interesting topic and deserves further study.

 \textbf{\emph{Proof of Theorem 2.2.} }
Note the fact from Lemma 1.8 of Stute (1993) that,
\[
F(t) - \bar F_n (t) =  - (1 - F(t))[\tilde \Lambda (t) - \Lambda (t)] + R'_{n1} (t) + R'_{n2} (t),
\]
where
\[
R'_{n1} (t) = 2^{ - 1} \exp \{  -\tilde \Lambda ^* (t)\} [\tilde \Lambda(t) - \Lambda (t)]^2,
\]
\[
R'_{n2} (t) = \exp \{  -\tilde \Lambda^{**} (t)\} [\tilde \Lambda(t) +\ln (1 - \bar F_n (t))]
\]
with $\tilde \Lambda^* (t)$ between $\tilde \Lambda(t)$ and  $\Lambda (t)$ and $\tilde \Lambda^{**} (t)$ between $\tilde \Lambda(t)$ and $ - \ln (1 - \bar F_n (t))$, respectively. Furthermore, following similar discussion to Lemma 3.2, we have
\[
\mathop {\sup }\limits_{0 \le x \le b} \left| {\tilde {\Lambda }(x) -
\Lambda (x)} \right|  = O(n^{ - 1 \mathord{\left/ {\vphantom {1 2}} \right.
\kern-\nulldelimiterspace} 2}(\log \log n)^{1 \mathord{\left/ {\vphantom {1
2}} \right. \kern-\nulldelimiterspace} 2})\ a.s.
\]
Hence, Lemma 3.4, 3.5 together with the result of Theorem 2.1 yield Theorem 2.2.  \ \ \ \ \ \ \   \qed

As an application of Theorem 2.2, one can obtain the LIL asymptotic result for the TWJ PL-estimator in the following.
In fact, some other similar results in Zhou and Yip(1999) can also be obtained by Theorem 2.1 or Theorem 2.2. We will consider these topics in future.

\textbf{\emph{Corollary 2.1.}} Suppose that $a_{F_A}  \le a_H$ and (A3) are satisfied. Then the stochastic sequence
\[\{ ({n \mathord{\left/
 {\vphantom {n {(2\log \log n)^{{1 \mathord{\left/
 {\vphantom {1 2}} \right.
 \kern-\nulldelimiterspace} 2}} (\tilde F_n (t) - F(t))}}} \right.
 \kern-\nulldelimiterspace} {(2\log \log n)^{{1 \mathord{\left/
 {\vphantom {1 2}} \right.
 \kern-\nulldelimiterspace} 2}} (\tilde F_n (t) - F(t))}})\}
 \]
  is almost surely relatively compact in the supermum norm of functions over $(a_H ,b]$, and its set of limit point is
\[
\{ (d(b))^{{1 \mathord{\left/
 {\vphantom {1 2}} \right.
 \kern-\nulldelimiterspace} 2}} (1 - F( \cdot ))g({{d( \cdot )} \mathord{\left/
 {\vphantom {{d( \cdot )} {d(b)}}} \right.
 \kern-\nulldelimiterspace} {d(b)}}):g \in G\}
\]
where $d(t) = \int_{a_H }^t {R^{-2} (u)dN(u)}$ and $G$ is Strassen's set of absolutely continuous functions,
\[
G = \{ g\left| g \right.:[0,1] \to R,\ g(0) = 0,\int_0^1 {(\frac{{dg(t)}}{{dt}})^2 dt}  \le 1\} .
\]
Consequently, write $v^2 (t) = (1 - F(t))d(t),$ then
 \[
\mathop {\lim \sup }\limits_{n \to \infty } (\frac{n}{{2\log \log n}})^{{1 \mathord{\left/
 {\vphantom {1 2}} \right.
 \kern-\nulldelimiterspace} 2}} \mathop {\sup }\limits_{a_H  < t \le b} \left| {\tilde F_n (t) - F(t)} \right| = \mathop {\sup }\limits_{a_H  < t \le b} v(t)\ a.s.,
\]
and
  \[
\mathop {\lim \inf }\limits_{n \to \infty } (n\log \log n)^{{1 \mathord{\left/
 {\vphantom {1 2}} \right.
 \kern-\nulldelimiterspace} 2}} \mathop {\sup }\limits_{a_H  < t \le b} \frac{{\left| {\tilde F_n (t) - F(t)} \right|}}{{1 - F(t)}} = \frac{\pi }{{2\sqrt 2 }}(d(b))^{{1 \mathord{\left/
 {\vphantom {1 2}} \right.
 \kern-\nulldelimiterspace} 2}}\ a.s.
\]
\textbf{\emph{Proof of Corollary 2.1.}} The rate of the strong convergence in Theorem 2.2 provides enough support to the result of Corollary 2.1, and the proof is similar to the procedure of Corollary 2.2 in Zhou and Yip (1999), we omit the details here.  \ \ \ \ \ \ \  \ \ \  \ \ \ \ \ \ \ \ \ \ \ \ \ \ \ \ \ \ \ \ \  \ \ \  \ \ \ \ \ \ \ \ \ \ \ \ \ \ \ \ \ \ \ \ \  \ \ \   \qed

\section{Some Lemmas and their proofs}
Let $M$ be a generic positive constant in the sequel, which could take different values at different places.

\textbf{\emph{Lemma 3.1.}}  Under the d.f. continuity of random variable assumed above, for $t<b_H$,
\[
\ \ \ \ \ \ \ \ \ \ \ \ \ \ \ \ \ \ \ \ \ \ \ \ \ \ \ \ \ \ \ \ \ \ \
\mathop {\sup }\limits_{0 \le x \le t} \left| {\int_0^x {(\tilde K^{ - 1}  - K^{ - 1} )d(\tilde Q_1  - Q_1 )} } \right| = O(n^{{{ - 3} \mathord{\left/
 {\vphantom {{ - 3} 4}} \right.
 \kern-\nulldelimiterspace} 4}} (\log n)^{{3 \mathord{\left/
 {\vphantom {3 4}} \right.
 \kern-\nulldelimiterspace} 4}} )\ a.s.,  \ \ \ \ \ \ \ \ \ \ \ \ \ \ \ \ \  (3.1)
\]
and
\[
 \ \ \ \ \ \ \ \ \ \ \ \ \ \ \ \ \ \ \ \ \ \ \ \ \ \ \ \ \ \ \ \ \ \ \
\mathop {\sup }\limits_{0 \le x \le t} \left| {\int_0^x {(\tilde K^{ - 1}  - K^{ - 1} )d(\tilde Q_2  - Q_2 )} } \right| = O(n^{{{ - 3} \mathord{\left/
 {\vphantom {{ - 3} 4}} \right.
 \kern-\nulldelimiterspace} 4}} (\log n)^{{3 \mathord{\left/
 {\vphantom {3 4}} \right.
 \kern-\nulldelimiterspace} 4}} )\ a.s. \ \ \ \ \ \ \ \ \ \ \ \ \ \ \ \ \ \ \ (3.2)
\]

\textbf{\emph{Proof of Lemma 3.1.}} Partitioning firstly the interval $[0,x]$ into subintervals $[x_i ,x_{i + 1} ],i = 1,2,
\cdots ,k_n $, with $k_n = O({{\sqrt n } \mathord{\left/
 {\vphantom {{\sqrt n } {(\log n)^{{1 \mathord{\left/
 {\vphantom {1 2}} \right.
 \kern-\nulldelimiterspace} 2}} }}} \right.
 \kern-\nulldelimiterspace} {(\log n)^{{1 \mathord{\left/
 {\vphantom {1 2}} \right.
 \kern-\nulldelimiterspace} 2}} }})$,
and $0 = x_1 < x_2 < \cdots < x_{k_n + 1} = x$ such that $\max_{1\le i\le k_n}|x_{i+1}-x_{i}|\le O ((n/\log n)^{-1/2})$. Because   $K_1(\cdot)$ is continuous differentable function, then
\[
K_1 (x_i ) - K_1 (x_{i + 1} )  =
 O({{(\log n)^{{1 \mathord{\left/
 {\vphantom {1 2}} \right.
 \kern-\nulldelimiterspace} 2}} } \mathord{\left/
 {\vphantom {{(\log n)^{{1 \mathord{\left/
 {\vphantom {1 2}} \right.
 \kern-\nulldelimiterspace} 2}} } {\sqrt n }}} \right.
 \kern-\nulldelimiterspace} {\sqrt n }}) , \
K_2 (x_i ) - K_2 (x_{i + 1} )  =  O({{(\log n)^{{1 \mathord{\left/
 {\vphantom {1 2}} \right.
 \kern-\nulldelimiterspace} 2}} } \mathord{\left/
 {\vphantom {{(\log n)^{{1 \mathord{\left/
 {\vphantom {1 2}} \right.
 \kern-\nulldelimiterspace} 2}} } {\sqrt n }}} \right.
 \kern-\nulldelimiterspace} {\sqrt n }}).
\]
Then $
K(x_i ) - K(x_{i + 1} )  = O({{(\log n)^{{1 \mathord{\left/
 {\vphantom {1 2}} \right.
 \kern-\nulldelimiterspace} 2}} } \mathord{\left/
 {\vphantom {{(\log n)^{{1 \mathord{\left/
 {\vphantom {1 2}} \right.
 \kern-\nulldelimiterspace} 2}} } {\sqrt n }}} \right.
 \kern-\nulldelimiterspace} {\sqrt n }}).$

Note that  $\tilde {Q}_1 (x)$ and $Q_1 (x)$ are monotone increasing function, we have, as in the proof of Lemma 2 of Lo and Singh (1986), that the left hand side in (3.1) is bounded  by
\begin{eqnarray*}
& &  \left| {\int_0^x {(\tilde K^{ - 1}  - K^{ - 1} )d(\tilde Q_1  - Q_1 )} } \right| \\
  & \le &  \mathop {\sup }\limits_{y \in [0,x]} \left| {\tilde K^{ - 1} (y) - K^{ - 1} (y)} \right|\sum\limits_{i = 1}^{k_n } {\left| {[\tilde Q_1 (x_{i + 1} ) - Q_1 (x_{i + 1} )] - [\tilde Q_1 (x_i ) - Q_1 (x_i )]} \right|}  \\
 & &  + \sum\limits_{i = 1}^{k_n } {\int_{x_i }^{x_{i + 1} } {\left| {[\tilde K^{ - 1} (x) - K^{ - 1} (x)] - [\tilde K^{ - 1} (x_i ) - K^{ - 1} (x_i )]} \right|d \tilde Q_1 (x)} }\\
 & &  + \sum\limits_{i = 1}^{k_n } {\int_{x_i }^{x_{i + 1} } {\left| {[\tilde K^{ - 1} (x) - K^{ - 1} (x)] - [\tilde K^{ - 1} (x_i ) - K^{ - 1} (x_i )]} \right|d Q_1 (x)} }\\
  &  \le & 2\mathop {\max }\limits_{1 \le i \le k_n } \mathop {\sup }\limits_{y \in [x_i ,x_{i + 1} ]} \left| {[\tilde K^{ - 1} (x_i ) - K^{ - 1} (x_i )] - [\tilde K^{ - 1} (y) - K^{ - 1} (y)]} \right| \\
   & &  + k_n \mathop {\max }\limits_{1 \le i \le k_n } \left| {\tilde Q_1 (x_{i + 1} ) - Q_1 (x_{i + 1} ) - (\tilde Q_1 (x_i ) - Q_1 (x_i ))} \right|\mathop {\sup }\limits_{y \in [0,x]} \left| {\tilde K^{ - 1} (y) - K^{ - 1} (y)} \right| \\
 &  = :&A + B.
\end{eqnarray*}
For estimating A, we further subdivide every $[x_i ,x_{i + 1} ]$ into
subintervals $[x_{ij} ,x_{i(j + 1)} ], j = 1, \cdots ,a_n ,$ with $a_n  = O(n^{{1 \mathord{\left/
 {\vphantom {1 4}} \right.
 \kern-\nulldelimiterspace} 4}} \log ^{{{ - 1} \mathord{\left/
 {\vphantom {{ - 1} 4}} \right.
 \kern-\nulldelimiterspace} 4}} n)$, such that
\[
\begin{array}{l}
 K_1 (x_{ij} ) - K_1 (x_{i(j + 1)} ) =  O(n^{{{ - 3} \mathord{\left/
 {\vphantom {{ - 3} 4}} \right.
 \kern-\nulldelimiterspace} 4}} (\log n)^{{3 \mathord{\left/
 {\vphantom {3 4}} \right.
 \kern-\nulldelimiterspace} 4}} );\begin{array}{*{20}c}\\
\end{array}K_2 (x_{ij} ) - K_2 (x_{i(j + 1)} ) =  O(n^{{{ - 3} \mathord{\left/
 {\vphantom {{ - 3} 4}} \right.
 \kern-\nulldelimiterspace} 4}} (\log n)^{{3 \mathord{\left/
 {\vphantom {3 4}} \right.
 \kern-\nulldelimiterspace} 4}} ),
  \end{array}
\]
uniformly in $i,j$. Now, since $\sup \left| {\tilde K - K} \right|^2  = O(n^{ - 1} \log n )\ a.s.
$ by LIL, and  $K(\cdot)$ and $\tilde {K}(\cdot)$ are bound in the intervals $[x_i ,x_{i + 1} ], i = 1,2,
\cdots ,k_n, $ it follows that
\begin{eqnarray*}
 & & \mathop {\sup }\limits_{y \in [x_i ,x_{i + 1} ]} \left| {[\tilde K^{ - 1} (x_i ) - K^{ - 1} (x_i )] - [\tilde K^{ - 1} (y) - K^{ - 1} (y)]} \right| \\
     &  \le &  \mathop {\sup }\limits_{y \in [x_i ,x_{i + 1} ]} \{ K^{ - 2} (x_{i + 1} )\left| {[\tilde K(y) - K(y)] - [\tilde K(x_i ) - K(x_i )]} \right| \\
 & & + \left| {\frac{{K^2 (x_{i + 1} ) - K^2 (y)}}{{K^2 (y)K^2 (x_{i + 1} )}}} \right|\left| {\tilde K(y) - K(y)} \right| \\
 & &  + \left| {\frac{{K^2 (x_{i + 1} ) - K^2 (x_i )}}{{K^2 (x_i )K^2 (x_{i + 1} )}}} \right|\left| {\tilde K(x_i ) - K(x_i )} \right| + O(n^{ - 1} \log n )\} \ a.s.\\
  & \le & M\mathop {\max }\limits_{1 \le j \le a_n } \left| {\tilde K(x_{ij} ) - \tilde K(x_i ) - K(x_{ij} ) + K(x_i )} \right| + O(n^{{{ - 3} \mathord{\left/
 {\vphantom {{ - 3} 4}} \right.
 \kern-\nulldelimiterspace} 4}} (\log n)^{{3 \mathord{\left/
 {\vphantom {3 4}} \right.
 \kern-\nulldelimiterspace} 4}} )\ a.s.
\end{eqnarray*}

Set $\eta _k = I(x_i < a_k \le x_{ij} ) - P(x_i < A \le x_{ij} )$, if we take $c=1$, $\sigma ^2  = M
n^{{{ - 1} \mathord{\left/
 {\vphantom {{ - 1} 2}} \right.
 \kern-\nulldelimiterspace} 2}} (\log n)^{{1 \mathord{\left/
 {\vphantom {1 2}} \right.
 \kern-\nulldelimiterspace} 2}},z = M\log n$, then the following probability bound can be verified from the exponential inequality of Lemma 1 in Lo and Singh (1986).
 \[
\mathop {\max }\limits_{i \le k_n } \mathop {\max }\limits_{1 \le j \le a_n } P(\left| {\tilde K(x_{ij} ) - \tilde K(x_i ) - K(x_{ij} ) + K(x_i )} \right| > Mn^{{{ - 3} \mathord{\left/
 {\vphantom {{ - 3} 4}} \right.
 \kern-\nulldelimiterspace} 4}} (\log n)^{{3 \mathord{\left/
 {\vphantom {3 4}} \right.
 \kern-\nulldelimiterspace} 4}} ) = O(n^{ - 3} ).
\]

Utilizing Bonferroni inequality together with the Borel Cantelli Lemma, it follows that $A = O(n^{{{ - 3} \mathord{\left/
 {\vphantom {{ - 3} 4}} \right.
 \kern-\nulldelimiterspace} 4}} (\log n)^{{3 \mathord{\left/
 {\vphantom {3 4}} \right.
 \kern-\nulldelimiterspace} 4}} )
\ a.s.$ The estimation of $B$ is treated similarly and leads to the same order. The proof of (3.2) is similar, we omit the details here. This completes the proof.  \ \  \ \ \   \qed

Define $ \tilde \Lambda _A (t)  = \int_0^t {\tilde K^{-1}(u)d\tilde Q(u)} $ as the nonparametric estimate for $\Lambda _A (\cdot)$, the cumulative hazard function of $A$.

\emph{\textbf{Lemma 3.2.}} When $b < b_H $, we have
\begin{eqnarray*}
&& \mathop {\sup }\limits_{0 \le x \le b} \left| {\tilde {\Lambda }_A (x) -
\Lambda _A (x)} \right| = O(n^{ - 1 \mathord{\left/ {\vphantom {1 2}}
\right. \kern-\nulldelimiterspace} 2}(\log \log n)^{1 \mathord{\left/
{\vphantom {1 2}} \right. \kern-\nulldelimiterspace} 2})\ a.s.
\end{eqnarray*}

\emph{\textbf{Proof of Lemma 3.2.}} Denote
\begin{eqnarray*}
I_{n1} &=& \mathop {\sup }\limits_{0 \le x \le b} \left| {\int_0^x {\frac{d\tilde
{Q}(u)}{\tilde {K}(u)}} - \int_0^x {\frac{d\tilde {Q}(u)}{K(u)}} } \right|, \  I_{n2} =
\mathop {\sup }\limits_{0 \le x \le b} \left| {\int_0^x {\frac{d\tilde
{Q}(u)}{K(u)}} - \int_0^x {\frac{dQ(u)}{K(u)}} } \right|.
 \end{eqnarray*}

As to $I_{n1}$, combining the condition  $b < b_H $ and LIL for empirical processes, there is
  \begin{eqnarray*}
 I_{n1}  & \le& \mathop {\sup }\limits_{0 \le x \le b} (\mathop {\sup }\limits_{0 \le u
\le x} [\tilde {K}(u) - K(u)]^2)\int_0^x {\frac{d\tilde {Q}(u)}{K^2(u)\tilde
{K}(u)}} \\
 &&+ \mathop {\sup }\limits_{0 \le x \le b} \int_0^x {\left| {\frac{\tilde
{K}(u) - K(u)}{K^2(u)}} \right|d\tilde {Q}(u)} = O((n^{ - 1} \log \log n)^{{1 \mathord{\left/
 {\vphantom {1 2}} \right.
 \kern-\nulldelimiterspace} 2}} ) \ a.s.
 \end{eqnarray*}
Next, put
\[
 J_n (x) = \int_0^x K^{-1}(u) d[\tilde {Q}(u) - Q(u)],
\]
then $I_{n2} = \mathop {\sup }\limits_{0 \le x \le b} J_n (x)$ is of the order $(n^{-1}{\log \log n})^{1 \mathord{\left/
{\vphantom {1 2}} \right. \kern-\nulldelimiterspace} 2}$ almost sure.

 In fact, the process $J_n (x)$ satisfies LIL, 
since it is an empirical process over VC classes of function with square integral envelope, and thus $\mathop {\sup }\limits_{0 < x \le b} \left| {\tilde {\Lambda }_A (x) - \Lambda _A (x)} \right|$ is also of the same order. This ends the proof. \ \ \ \ \ \ \ \ \ \ \ \ \ \ \  \ \ \  \ \ \ \ \ \ \ \ \ \ \ \ \ \ \  \ \ \ \ \ \ \ \ \ \ \ \ \ \ \ \ \ \  \ \ \ \ \ \ \ \ \ \ \ \ \ \ \ \ \ \  \ \ \ \ \ \ \ \ \ \ \ \ \ \ \ \ \ \  \ \ \  \qed

We now establish a strong representation for $\tilde {S}_A(\cdot) ,$
which is constructed by pooling data from the truncation time and the
observed residual survival time.

\textbf{\emph{Lemma 3.3}}. When $b < b_H $, the stochastic process $\tilde {S}_A (t) - S_A (t)$ has an
asymptotic representation
\[
\tilde {S}_A (t) - S_A (t) =  n^{ - 1}\sum\limits_{i = 1}^n {S_A (t)\phi _i (t)}+
R_{n3} (t),
\]
 where $ \mathop {\sup }\limits_{a_H  \le t \le b} \left| {R_{n3} (t)} \right| = O(n^{{{ - 3} \mathord{\left/
 {\vphantom {{ - 3} 4}} \right.
 \kern-\nulldelimiterspace} 4}} (\log n)^{{3 \mathord{\left/
 {\vphantom {3 4}} \right.
 \kern-\nulldelimiterspace} 4}} )
\ a.s.$

\emph{\textbf{Remark 3}}  Lemma 3.3 indicates that $\tilde {S}_A ( \cdot )$ is a strong consistent estimator of $S_A ( \cdot )$, and obviously it implies the asymptotic representation for $\tilde {R}( \cdot ),$
\begin{eqnarray*}
 \tilde R(t) &= & \bar R(t) + \frac{1}{n}\sum\limits_{i = 1}^n {\{ I(a_i  > t) - S_A (t) - S_A (t)\phi _i (t)\} }  +O(n^{{{ - 3} \mathord{\left/
 {\vphantom {{ - 3} 4}} \right.\kern-\nulldelimiterspace} 4}} (\log n)^{{3 \mathord{\left/  {\vphantom {3 4}} \right. \kern-\nulldelimiterspace} 4}} )\ a.s.
 \end{eqnarray*}
where $\bar R(t) = n^{ - 1}\sum\nolimits_{i = 1}^n {I(a_i \le t \le y_i )} $.

\textbf{\emph{Proof of Lemma 3.3.}}
 By the definition of $\phi _i (t)$, there is
\[
- n^{ - 1} \sum\limits_{i = 1}^n {\phi _i (t)}  =  - \int_0^t {\frac{{\tilde K(u)}}{{K^2 (u)}}} dQ(u) + \int_0^t {\frac{d\tilde Q(u)}{{K(u)}}}.
\]

Since again
\begin{eqnarray*}
&& \left| {\int_0^t {(\frac{1}{\tilde {K}(u)}} - \frac{1}{K(u)})d(\tilde
{Q}(u) - Q(u))} \right| \\
  & \le&  \left| {\int_0^t {(\frac{1}{\tilde {K}(u)}} - \frac{1}{K(u)})d(\tilde
{Q}_1 (u) - Q_1 (u))} \right| + \left| { \int_0^t {(\frac{1}{\tilde {K}(u)}}-
\frac{1}{K(u)})d(\tilde {Q}_2 (u) - Q_2 (u))} \right|\\
& = :& S_1  + S_2.
 \end{eqnarray*}

Thus applying Lemma 3.1 to $S_1$  and $ S_2$, one can derive the following asymptotic representation under $b < b_H $,
\begin{eqnarray*}
 && \tilde \Lambda _A (t) - \Lambda _A (t)= \int_0^t {\frac{{d\tilde Q(u)}}{{\tilde K(u)}}}  - \int_0^t {\frac{{dQ(u)}}{{K(u)}}}\\
         &\le&  - \frac{1}{n}\sum\limits_{i = 1}^n {\phi _i (t)}  + \sup_{0\le u\le b} (K(u) - \tilde K(u))^2 \int_0^t {\frac{dQ(u)}{{\tilde K(u)K^2 (u)}}}  +O(n^{{{ - 3} \mathord{\left/
 {\vphantom {{ - 3} 4}} \right.
 \kern-\nulldelimiterspace} 4}} (\log n)^{{3 \mathord{\left/
 {\vphantom {3 4}} \right.
 \kern-\nulldelimiterspace} 4}} )\ a.s.
 \\
   &=& - \frac{1}{n}\sum\limits_{i = 1}^n {\phi _i (t)}+O(n^{{{ - 3} \mathord{\left/
 {\vphantom {{ - 3} 4}} \right.
 \kern-\nulldelimiterspace} 4}} (\log n)^{{3 \mathord{\left/
 {\vphantom {3 4}} \right.
 \kern-\nulldelimiterspace} 4}} )
\ a.s.,
\end{eqnarray*}
where
\begin{eqnarray*}
&& \sup_{0\le u\le b} (K(u) - \tilde K(u))^2 \left| {\int_0^t {\frac{dQ(u)}{{\tilde K(u)K^2 (u)}}} } \right| \\
&=& O(n^{ - 1} \log n)[ \left| {\int_0^t {\frac{dQ(u)}{{K^3 (u)}}} } \right| + \sup_{0\le u\le b} \left| {\frac{{K(u) - \tilde K(u)}}{{\tilde K(u)}}} \right|\left| {\int_0^t {\frac{dQ(u)}{{K^3 (u)}}} } \right|]\ a.s. \\
  &=&  O(n^{ - 1} \log n)\ a.s.
\end{eqnarray*}

Using Lemma 3.2, by expansion of the function exp$\{-x\}$ in neighborhood of zero.
 \begin{eqnarray*}
  && \tilde S_A (t) - S_A (t) = \exp \{  - \tilde \Lambda _A (t)\}  - \exp \{  - \Lambda _A (t)\}  \\
  &=&-\exp \{  - \Lambda _A (t)\} [ - \frac{1}{n}\sum\limits_{i = 1}^n {\phi _i (t)}  + O(n^{{{ - 3} \mathord{\left/
 {\vphantom {{ - 3} 4}} \right.
 \kern-\nulldelimiterspace} 4}} (\log n)^{{3 \mathord{\left/
 {\vphantom {3 4}} \right.
 \kern-\nulldelimiterspace} 4}} )
 +  O(n^{ - 1} \log \log n)]\ a.s.\\
 & =& \frac{1}{n}\sum\limits_{i = 1}^n {S_A (t)\phi _i (t)}  + O(n^{{{ - 3} \mathord{\left/
 {\vphantom {{ - 3} 4}} \right.
 \kern-\nulldelimiterspace} 4}} (\log n)^{{3 \mathord{\left/
 {\vphantom {3 4}} \right.
 \kern-\nulldelimiterspace} 4}} )
\ a.s.
 \end{eqnarray*}
 This ends the proof. \ \ \ \ \ \ \ \ \ \ \ \ \ \ \ \ \ \ \ \ \  \ \ \  \ \ \ \ \ \ \ \ \ \ \ \ \ \ \ \ \ \ \ \ \  \ \ \  \ \ \ \ \ \ \ \ \ \ \ \ \ \ \ \ \ \ \ \ \  \ \ \  \ \ \ \ \ \ \ \ \ \ \ \ \ \ \ \ \ \ \ \ \  \ \ \  \ \ \ \qed

 Next, similar to the discussion of some lemmas in Zhou and Yip (1999), we may derive two relevant lemmas under LBRC mechanism. Note that for $0 < b < b_H $, it follows from the SLLN
that
\[
\ \ \ \ \ \ \ \ \ \ \ \ \ \ \ \ \ \ \ \ \ \ \ \ \ \ \ \ \ \ \ \ \ \ \ \ \ \ \ \ \int_0^b {\frac{d\bar {N}(u)}{R(u)[\tilde {R}(u) + n^{ - 1}]}} < \infty .\ \ \ \ \ \ \ \ \ \ \ \ \ \ \ \ \ \ \ \ \ \ \ \ \ \ \ \ \ \ \ \ \ \ \ \ \ \ (3.3)
\]

For the proof of Theorem 2.2, we need  a slight modification of $\tilde {F}_n( \cdot) $. Define a new estimator $\bar {F}_n( \cdot) $ as
\[
1 - \bar {F}_n (x)  = \prod\limits_{y_i \le x} {[1 -
\frac{1}{n\tilde {R}(y_i ) + 1}]^{\delta _i }},
\]
which is only to safeguard against log0 when taking logarithms of $1 - \tilde {F}_n (x)$.

\textbf{\emph{Lemma 3.4.}} Under (A3), when $b < b_H $, there is
\[
\mathop {\sup }\limits_{0 \le x \le b} \left| {\tilde {F}_n(x)-\bar {F}_n (x)} \right| = O(n^{ - 1}) \ a.s.
\]

\textbf{\emph{Proof of Lemma 3.4.}} Obviously, with (3.3),
\begin{eqnarray*}
&& \mathop {\sup }\limits_{0 \le x \le b} \left| {\tilde
{F}_n (x)}- \bar {F}_n (x) \right| \\
    & \le& \frac{1}{n}\mathop {\sup }\limits_{0 \le x \le b} \int_0^x {\frac{{d\bar N(u)}}{{R(u)[\tilde R(u) + n^{ - 1} ]}}}  + \frac{1}{n}\mathop {\sup }\limits_{0 \le x \le b} \int_0^x {\left| {\frac{{R(u) - \tilde R(u)}}{{R(u)\tilde R(u)[\tilde R(u) + n^{ - 1} ]}}} \right|d\bar N(u)}  \\
    &\le& \frac{1}{n}\int_0^b {\frac{{d\bar N(u)}}{{R(u)[\tilde R(u) + n^{ - 1} ]}}}  + \frac{1}{n}\sup_{0\le u\le b} \frac{{\left| {R(u) - \tilde R(u)} \right|}}{{\tilde R(u)}}\int_0^b {\frac{{d\bar N(u)}}{{R(u)[\tilde R(u) + n^{ - 1} ]}}}  \\
  &=& O(n^{ - 1} )\ a.s.
  \end{eqnarray*}

This completes the proof. \ \ \ \ \ \ \ \ \ \ \ \ \ \ \ \ \ \ \ \ \  \ \ \  \ \ \ \ \ \ \ \ \ \ \ \ \ \ \ \ \ \ \ \ \  \ \ \  \ \ \ \ \ \ \ \ \ \ \ \ \ \ \ \ \ \ \ \ \  \ \ \ \ \ \ \ \ \ \ \ \  \ \ \   \qed

\textbf{\emph{Lemma 3.5.}} Under (A3), when $b < b_H $, there is
\[
\mathop {\sup }\limits_{0 \le x \le b} \left| {\tilde {\Lambda }(x) +\log (1 - \bar {F}_n (x))
} \right| = O(n^{ - 1})\ a.s.
\]

\textbf{\emph{Proof of Lemma 3.5.}} Similar to the discussion in Lemma 3.4, and using the Taylor's expansion for the function $\log (1 - x) $ when $x<1$, we have
\begin{eqnarray*}
 \mathop {\sup }\limits_{0 \le x \le b} \left| {\tilde {\Lambda }(x) +\log (1 - \bar {F}_n (x))
} \right|  & = &\mathop {\sup }\limits_{0 \le x \le b} \left| {\sum\limits_{i =
1}^n {\frac{\delta _i I(y_i \le x)}{n\tilde {R}(y_i )}} } + \sum\limits_{i:y_i \le x}
{\delta _i \log [1 - \frac{1}{n\tilde {R}(y_i ) + 1}]}\right| \\
 &\le& \mathop {\sup }\limits_{0 \le x \le b} \sum\limits_{i:y_i \le x}
{\frac{\delta _i}{n\tilde {R}(y_i )[n\tilde {R}(y_i ) + 1]}} = O(n^{
- 1})\ a.s.
  \end{eqnarray*}

This completes the proof. \ \ \ \ \ \ \ \ \ \ \ \ \ \ \ \ \ \ \ \ \  \ \ \  \ \ \ \ \ \ \ \ \ \ \ \ \ \ \ \ \ \ \ \ \  \ \ \  \ \ \ \ \ \ \ \ \ \ \ \ \ \ \ \ \ \ \ \ \  \ \ \ \ \ \ \ \ \ \ \ \  \ \ \   \qed

\small
\baselineskip=8pt

\end{document}